\documentclass[preprint, 12pt]{elsarticle}

\usepackage{amsmath, amssymb, amsfonts}
\usepackage{amsthm}
\usepackage{mathtools}

\usepackage[T1]{fontenc}
\usepackage[utf8]{inputenc}
\usepackage{lmodern}

\usepackage{graphicx}
\usepackage{booktabs}
\usepackage{xcolor}

\theoremstyle{definition}
\newtheorem{definition}{Definition}[section]
\newtheorem{lemma}[definition]{Lemma}
\newtheorem{proposition}[definition]{Proposition}

\theoremstyle{plain}

\theoremstyle{remark}

\usepackage[colorlinks=true, linkcolor=black, citecolor=black, urlcolor=black]{hyperref}

\begin{document}
\begin{frontmatter}

\title{Spectral Reconstruction in Fractional Derivative Order}

\author{Derek C. Braun}
\address{School of Science, Technology, Accessibility, Mathematics and Public Health, Gallaudet University, Washington, DC 20002, USA}
\ead{derek.braun@gallaudet.edu}

\begin{abstract}
We introduce the spectral reconstruction operator $\mathcal{S}$, which treats derivative order as a continuous variable and integrates fractional derivative data across all nonnegative real orders to represent analytic functions. Because the fractional derivatives in the integrand draw on function values across an infinite stretch of the real line, the spectral integrand encodes information about $f$ that no individual integer derivative at the origin can capture. We define an admissible class of fractional derivative families that recover the classical derivative ladder, satisfy the semigroup property, and keep $D^r f(0)$ finite (a condition that excludes monomials). We prove that the spectral reconstruction operator asymptotically reconstructs the exponential function as $x\to\infty$. We further show that the spectral operator and the classical Maclaurin series are formed from the same underlying function of derivative order: the spectral operator integrates it continuously, while the Maclaurin series samples it at the nonnegative integers. Classical summation formulas therefore quantify this discrepancy under their respective hypotheses. Under their respective hypotheses, the Euler--Maclaurin formula gives an exact identity consisting of a finite correction series and its remainder, while the Abel--Plana formula gives an exact correction integral. Examples drawn from entire, oscillatory, rapidly decaying, special-function, and finite-radius analytic classes demonstrate stable reconstruction within their applicable domains. The Abel--Plana representation matches $f$ to the full precision of our computations in every case. In the finite-radius examples studied here, one with a real pole and one with a complex-conjugate pole pair, the spectral integral and the Maclaurin series encounter the same convergence boundary, although through different mathematics. In the examples studied, the sum--integral discrepancy remains bounded, so the uncorrected operator provides useful approximations where $f(x)$ is large. The leading summation term $\mathcal{E}_0=\frac{1}{2}f(0)$ and higher terms then systematically account for the discrepancy concentrated near the origin. These results establish derivative order as a meaningful, non-local continuous domain for the representation and reconstruction of analytic functions.
\end{abstract}

\end{frontmatter}

\section{Introduction}
Fractional calculus extends differentiation to continuous order. If the derivative of a function can be taken to any real order, a natural question follows: can a function be reconstructed from its fractional derivative data integrated continuously over real orders, rather than sampled only at the integers?

The idea of fractional derivatives dates back to a 1695 letter from Leibniz to l'Hôpital, asking what a half-order derivative would mean \cite{Ross1977}. Classical constructions extend differentiation to non-integer order in several ways, including the Riemann--Liouville and Liouville (Fourier-multiplier) \cite{Weyl1917} definitions. Osler proved a generalized Taylor theorem reconstructing $f$ from discrete fractional derivative data \cite{Osler1971}, a striking result showing that fractional order carries deep analytic information. Subsequent work has explored related fractional Taylor and Maclaurin series under a variety of definitions and constructions \cite{Alquran2023, FernandezBaleanu2018, Gladkina2017, Wei2020, ZitaneTorres2023}. All of these constructions are discrete in derivative order. To our knowledge, no prior work has asked and answered whether a continuous integral over derivative order can reconstruct analytic functions.

We introduce the \textit{spectral reconstruction operator}
\begin{equation}
\mathcal{S}[f](x)
:=
\int_0^\infty
\frac{D^r f(0)\,x^r}{\Gamma(r+1)}
\,dr,
\end{equation}
which integrates fractional derivative data across all nonnegative real orders. Because the fractional derivatives in the integrand draw on function values across an infinite stretch of the real line, the spectral operator encodes information about $f$ that no individual integer derivative at the origin can capture. Just as Fourier analysis treats frequency as a continuous variable, $\mathcal{S}$ treats derivative order as a continuous variable and reconstructs $f$ from its derivative-order spectrum.

We define an admissible class of fractional derivative families and prove that
\[
\mathcal S[e^x](x)\sim e^x \qquad (x\to\infty),
\]
the first proved instance of asymptotic reconstruction by the operator. A sum--integral identity relates the spectral operator to the classical Maclaurin series: the Maclaurin summand and the spectral integrand are both evaluations of the same underlying function of derivative order, one sampled at integers and one integrated continuously. Classical summation formulas are built to measure exactly this kind of discrepancy, under their respective hypotheses. We show that the Euler--Maclaurin summation formula gives an exact identity between $\mathcal{S}[f]$ and the Maclaurin series. The Abel--Plana formula gives an exact integral representation of the same discrepancy under stronger analyticity conditions.

Examples from entire, oscillatory, rapidly decaying, special-function, and finite-radius analytic classes confirm stable reconstruction within each function's applicable domain. In particular, the Abel--Plana representation matches $f$ to the full precision of our computations in every case. This establishes derivative order as a genuinely non-local, continuous domain for representing and reconstructing analytic functions.

\section{The Spectral Reconstruction Operator}
\label{sec:operator}
We now give this construction its formal definition: an integral over fractional derivative order that replaces the discrete sum of the classical Maclaurin series.

\begin{definition}
Let $D^{r}f(0)$ denote an admissible extension of the integer-order derivative data to real order $r\geq 0$. For real $x>0$ such that the integral converges, the \textit{spectral reconstruction operator} is
\begin{equation}
\label{eq:operator_def}
\mathcal{S}[f](x)
:=
\int_{0}^{\infty}
\frac{D^{r}f(0)\,x^{r}}{\Gamma(r+1)}
\,dr.
\end{equation}
The operator therefore depends on the chosen non-integer-order extension, which the admissibility conditions of Section~\ref{sub:frac_deriv_admissibility} specify.
\end{definition}

\subsection{Admissibility Conditions on the Fractional Derivative}
\label{sub:frac_deriv_admissibility}
We call a family $\{D^r\}_{r\ge0}$ \emph{admissible} if it satisfies the following conditions. We do not claim these conditions are minimal or necessary. They define the class considered here and formalize the interpretation of $r$ as accumulated differentiation order.

First, the family must recover the classical derivative ladder at integer orders:
\[
D^n f(0)=f^{(n)}(0),
\qquad n\in\mathbb N_0.
\]
Second, it must represent accumulated differentiation, satisfying the semigroup property:
\[
D^\alpha D^\beta = D^{\alpha+\beta}.
\]
Third, $D^r f(0)$ must exist and remain finite for real $r\geq 0$.

These requirements are restrictive: which ones hold depends on both the derivative definition and the function under consideration. In the examples below, we use the left-sided Liouville derivative based at $-\infty$, and the Fourier-multiplier constructions on $\mathbb R$,
\begin{equation}
\label{eq:liouville_fourier}
\mathcal F[D^r f](\xi)
=
(i\xi)^r\widehat f(\xi),
\end{equation}
with a fixed branch of the complex power. Other definitions are relevant only when they produce order data satisfying the same conditions.

Here is the heuristic behind the result: for large $x$, the integrand $x^r/\Gamma(r+1)$ is sharply peaked near derivative order $r=x$, so the reconstruction is dominated by orders comparable to $x$ itself.

\begin{proposition}[The Spectral Operator Meaningfully Approximates the Exponential]
\label{prop:meaningful_approximation_exponential}

Let $f(x)=e^x$, and suppose the chosen continuous-order derivative satisfies
\[
D^r e^x=e^x,
\qquad r\geq 0.
\]
Then
\[
\mathcal S[e^x](x)\sim e^x
\qquad \text{as } x\to\infty.
\]
\end{proposition}

\begin{proof}
We apply Stirling's formula to put the integral into the standard form for Laplace's method, locate this maximum precisely, and show that the resulting local Gaussian contribution has total mass asymptotic to $e^x$.

We use the left-sided Liouville derivative based at $-\infty$, which gives 
\[
\prescript{\mathrm{L}}{-\infty}{D}^{r} e^{x}= e^{x},
\qquad\text{so that}\qquad
\prescript{\mathrm{L}}{-\infty}{D}^{r} e^{x}\big|_{x=0} = 1,
\]
and
\[
\mathcal S[e^x](x)
=
\int_0^\infty \frac{x^r}{\Gamma(r+1)}\,dr.
\]
For each fixed $x>0$, the integral converges by Stirling's formula,
\[
\Gamma(r+1)
\sim
\sqrt{2\pi r}\left(\frac{r}{e}\right)^r.
\]
Hence
\[
\frac{x^r}{\Gamma(r+1)}
\sim
\frac{1}{\sqrt{2\pi r}}
\exp\!\bigl(\phi_x(r)\bigr),
\qquad
\phi_x(r)
=
r+r\log x-r\log r.
\]
The phase has a unique maximum at $r=x$, with
\[
\phi_x(x)=x,
\qquad
\phi_x''(x)=-\frac{1}{x}.
\]
Laplace's method therefore gives
\[
\mathcal S[e^x](x)
\sim
\frac{e^x}{\sqrt{2\pi x}}
\int_0^\infty
\exp\!\left(-\frac{(r-x)^2}{2x}\right)\,dr.
\]
As $x\to\infty$, the maximum lies a distance $\sqrt{x}$ standard deviations from the endpoint $r=0$, so the truncated Gaussian integral is asymptotic to the full Gaussian integral:
\[
\int_0^\infty
\exp\!\left(-\frac{(r-x)^2}{2x}\right)\,dr
\sim
\sqrt{2\pi x}.
\]
Therefore,
\[
\mathcal S[e^x](x)
\sim
\frac{e^x}{\sqrt{2\pi x}}
\cdot
\sqrt{2\pi x}
=
e^x.
\]
\end{proof}
\section{Sum--Integral Discrepancy}
\label{sec:sum_integral_discrepancy}
The spectral reconstruction operator $\mathcal{S}[f]$ integrates the integrand $\phi(r,x)$ continuously over all nonnegative real orders, while the Maclaurin series samples the same integrand at the nonnegative integers. Define the integrand
\[
\phi(r,x)
:=
\frac{D^r f(0)\,x^r}{\Gamma(r+1)}.
\]
\[
\mathcal S[f](x)
=
\int_0^\infty \phi(r,x)\,dr.
\]
Classical summation formulas quantify this discrepancy:
\[
\sum_{n=0}^{\infty} \phi(n,x)
=
\int_0^\infty \phi(r,x)\,dr
+
\text{correction terms}.
\]
At integer orders, the same function reproduces the summands of the Maclaurin series.

\begin{lemma}[Integer samples of the integrand are the Maclaurin summands]
\label{lem:integer_sampling}
For every integer $n\geq 0$,
\[
\phi(n,x)
=
\frac{f^{(n)}(0)x^n}{n!}.
\]
Consequently, whenever the Maclaurin series converges,
\[
\sum_{n=0}^\infty \phi(n,x)
=
\sum_{n=0}^\infty
\frac{f^{(n)}(0)x^n}{n!}.
\]
\end{lemma}

\begin{proof}
By admissibility,
\[
D^n f(0)=f^{(n)}(0),
\qquad n\in\mathbb N_0,
\]
and $\Gamma(n+1)=n!$. We substitute these into the definition of $\phi$ to obtain the result.
\end{proof}

Thus the spectral operator integrates $\phi(r,x)$ over continuous derivative order, whereas the Maclaurin series samples it at the nonnegative integers. Their discrepancy is therefore
\[
\sum_{n=0}^\infty \phi(n,x)
-
\mathcal S[f](x)
=
\sum_{n=0}^\infty \phi(n,x)
-
\int_0^\infty \phi(r,x)\,dr.
\]
Classical summation formulas analyze this sum--integral discrepancy under their respective regularity assumptions. We consider Euler--Maclaurin and Abel--Plana below, but this relationship is independent of any particular summation formula.

\subsection{Euler--Maclaurin summation formula}
The Euler--Maclaurin summation formula expresses the sum--integral discrepancy as a finite hierarchy of boundary corrections with an explicit remainder term~\cite[\S 2.10(i)]{NIST:DLMF}.

\begin{proposition}[The Euler--Maclaurin summation formula expresses the discrepancy through boundary corrections]
\label{prop:sum_integral_discrepancy}

Fix $x>0$ and $m\geq1$. Assume that $\phi(\cdot,x)$ and $\partial_r^{2m}\phi(\cdot,x)$ are integrable on $[0,\infty)$, that $\sum_{n=0}^{\infty}\phi(n,x)$ converges, and that $r\mapsto\phi(r,x)$ has $2m$ continuous derivatives on $[0,\infty)$, with
\[
\partial_r^j\phi(r,x)\longrightarrow 0
\qquad
(r\to\infty),
\qquad
0\leq j\leq 2m-1.
\]
Then
\[
\sum_{n=0}^{\infty}\phi(n,x)
=
\int_0^\infty \phi(r,x)\,dr
+
\frac12\phi(0,x)
-
\sum_{k=1}^{m}
\frac{B_{2k}}{(2k)!}
\frac{\partial^{2k-1}\phi}{\partial r^{2k-1}}(0,x)
+
R_m(x),
\]
where
\[
R_m(x)
=
-\frac{1}{(2m)!}
\int_0^\infty
\widetilde B_{2m}(r)
\frac{\partial^{2m}\phi}{\partial r^{2m}}(r,x)
\,dr,
\]
and
\[
\widetilde B_{2m}(r)
=
B_{2m}\bigl(r-\lfloor r\rfloor\bigr)
\]
is the periodic Bernoulli function. If, in addition, the Maclaurin series of $f$ converges to $f(x)$, then
\[
f(x)
=
\mathcal S[f](x)
+
\frac12 f(0)
-
\sum_{k=1}^{m}
\frac{B_{2k}}{(2k)!}
\frac{\partial^{2k-1}\phi}{\partial r^{2k-1}}(0,x)
+
R_m(x).
\]
\end{proposition}

\begin{proof}
We apply the Euler--Maclaurin formula to $r\mapsto\phi(r,x)$. The second identity then follows from Lemma~\ref{lem:integer_sampling} whenever the Maclaurin series converges to $f(x)$.
\end{proof}

The leading correction is
\[
\frac12\phi(0,x)=\frac12 f(0),
\]
and the higher terms depend on odd-order derivatives of the integrand at the boundary $r=0$. The identity is exact provided we retain $R_m(x)$; whether a finite truncation is useful depends on the size of this remainder. In the examples below, the first few correction terms reduce the discrepancy considerably.

\subsection{Abel--Plana formula}
The Abel--Plana formula requires stronger assumptions than those needed for the spectral reconstruction operator, but it gives an exact correction integral. In particular, it requires a complex-order continuation of the integrand,
\[
\phi(z,x)
=
\frac{D^z f(0)x^z}{\Gamma(z+1)},
\]
that is analytic in the right half-plane and satisfies the growth and integrability conditions required by the formula.

\begin{proposition}[Abel--Plana expresses the discrepancy as an exact correction integral]
\label{prop:abel_plana_discrepancy}

Fix $x>0$. Assume that $z\mapsto\phi(z,x)$ satisfies the hypotheses of the Abel--Plana formula. Then
\[
\sum_{n=0}^\infty \phi(n,x)
=
\int_0^\infty \phi(r,x)\,dr
+
\frac12\phi(0,x)
+
i
\int_0^\infty
\frac{\phi(it,x)-\phi(-it,x)}
{e^{2\pi t}-1}
\,dt.
\]
If, in addition, the Maclaurin series of $f$ converges to $f(x)$, then
\[
f(x)
=
\mathcal S[f](x)
+
\frac12 f(0)
+
i
\int_0^\infty
\frac{\phi(it,x)-\phi(-it,x)}
{e^{2\pi t}-1}
\,dt.
\]
\end{proposition}

\begin{proof}
We apply the Abel--Plana formula to $z\mapsto\phi(z,x)$. The second identity follows from Lemma~\ref{lem:integer_sampling} whenever the Maclaurin series converges to $f(x)$.
\end{proof}

Unlike the finite Euler--Maclaurin hierarchy, the Abel--Plana representation expresses the entire sum--integral discrepancy as a single correction integral. Both formulas contain the same leading boundary term,
\[
\frac12\phi(0,x)=\frac12 f(0).
\]
In the examples below, this term accounts for a substantial part of the observed discrepancy.

\begin{definition}[Correction notation]
\label{def:correction_series}

Define the leading correction by
\[
\mathcal E_0[f](x)
:=
\frac12 \phi(0,x)
=
\frac12 f(0).
\]
This term appears in both the Euler--Maclaurin and Abel--Plana formulas.

For $j\geq 1$, define the Euler--Maclaurin corrections by
\[
\mathcal E_j[f](x)
:=
-\frac{B_{2j}}{(2j)!}
\frac{\partial^{2j-1}\phi}{\partial r^{2j-1}}(0,x),
\]
where $B_{2j}$ denotes the Bernoulli numbers. The cumulative correction through index $N$ is
\[
\mathcal E_{\leq N}[f](x)
:=
\sum_{j=0}^{N}\mathcal E_j[f](x).
\]
\end{definition}

\section{Examples}
\label{sec:examples}
We evaluated the spectral reconstruction operator on representative analytic function classes to examine three questions: the behavior of the uncorrected operator, the effect of the common correction term $\mathcal{E}_0 = \frac12 f(0)$, and the extent to which higher-order Euler--Maclaurin terms or the Abel--Plana integral refine the reconstruction.

We performed all numerical evaluations using \texttt{mpmath} at 100-digit precision. We computed reconstruction integrals using adaptive quadrature (\texttt{mpmath.quad}).

The examples below span several classes of analytic functions:
\[
\begin{aligned}
f(x) &= e^x && \text{(exponential function)},\\
f(x) &= \sin x && \text{(oscillatory function)},\\
f(x) &= e^{-x^2} && \text{(rapidly decaying function)},\\
f(x) &= J_0(x) && \text{(special function with oscillatory decay)},\\
f(x) &= 1/(1-x) && \text{(function with a real pole at $x=1$ and $R=1$)},\\
f(x) &= 1/(1+x^2) && \text{(function with complex poles at $\pm i$ and $R=1$)}.
\end{aligned}
\]

In each example, we selected a fractional derivative from the admissible class of Section~\ref{sub:frac_deriv_admissibility}.

\subsection{Exponential function: $f(x)=e^x$}
\begin{figure}[tbp]
  \centering
  \includegraphics[width=\linewidth]{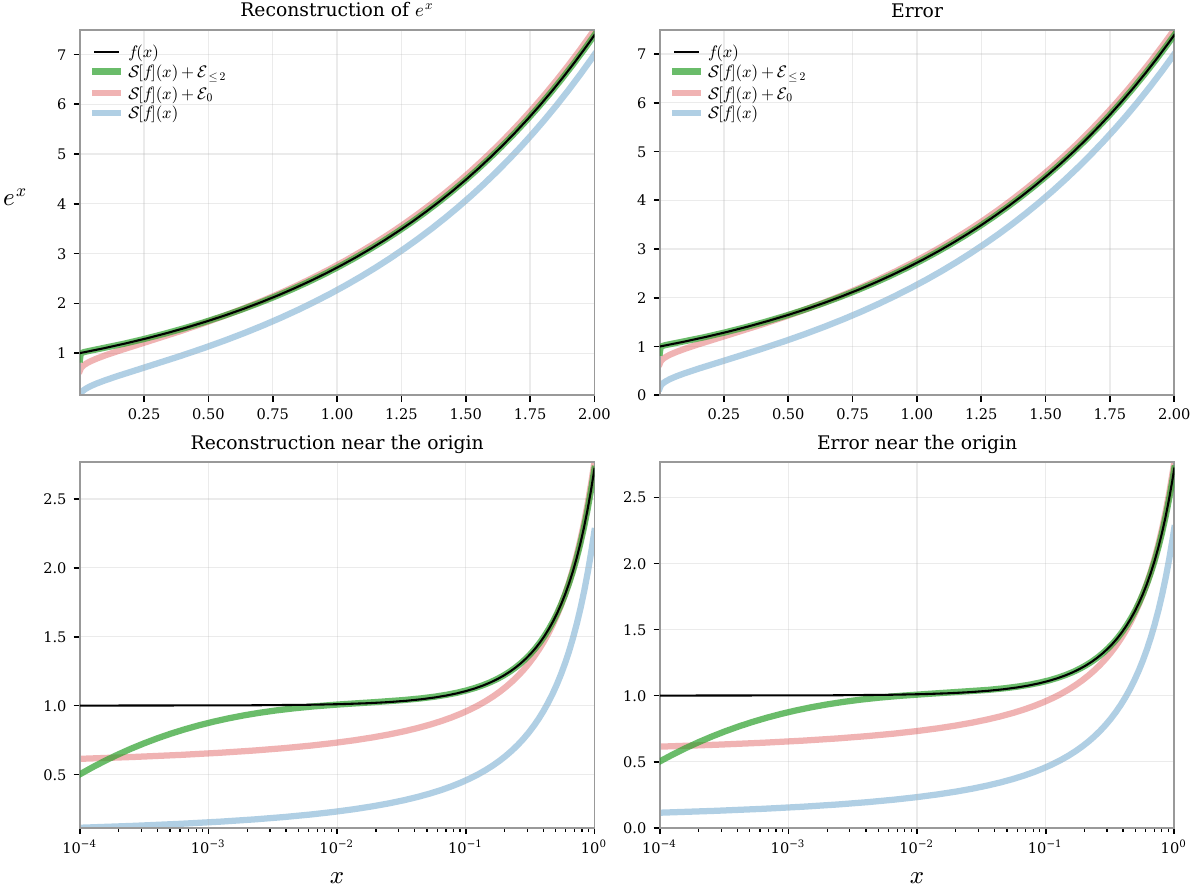}
  \caption{Spectral operator reconstruction of $e^x$ showing the uncorrected operator, with the common correction $\mathcal E_0$, and with Euler--Maclaurin corrections $\mathcal E_{\leq 2}$. The right panels show residuals.}
  \label{fig:exponential}
\end{figure}

\textbf{Fractional derivative}. For the exponential function, we use the left-sided Liouville derivative based at $-\infty$, which gives
\[
\prescript{\mathrm{L}}{-\infty}{D}^{r} e^{x}=e^{x},
\qquad\text{so that}\qquad
\prescript{\mathrm{L}}{-\infty}{D}^{r} e^{x}\big|_{x=0} = 1.
\]
\textbf{Correction terms}. The integrand simplifies to
\[
\phi(r,x) = \frac{x^{r}}{\Gamma(r+1)}.
\]
From Definition~\ref{def:correction_series}, the first three correction terms are:

\[
\mathcal{E}_{\leq 2}[e^{x}](x)
=
\underbracket{\frac{1}{2}}_{\textstyle \mathcal{E}_0}
\;+\;
\underbracket{- \frac{1}{12}\,(\log x + \gamma)}_{\textstyle \mathcal{E}_1}
\;+\;
\underbracket{\frac{1}{720}
     \!\left[(\log x+\gamma)^{3}
             - \frac{\pi^{2}}{2}(\log x+\gamma)
             + 2\zeta(3)\right]}_{\textstyle \mathcal{E}_2}
\;+\;
\cdots.
\]
Figure~\ref{fig:exponential} compares the uncorrected reconstruction, the common endpoint correction, and the fixed Euler--Maclaurin refinement $\mathcal E_{\leq 2}$.

\subsection{Oscillatory function: $f(x)=\sin x$}
\textbf{Fractional derivative}. For the sine function, the Liouville (Fourier-multiplier) definition is admissible and yields the phase-shift result below.
\[
{}^{\mathrm{L}}_{\mathbb{R}}D^{r} \sin x
= \sin\!\left(x + \frac{\pi r}{2}\right),
\qquad\text{so that}\qquad
{}^{\mathrm{L}}_{\mathbb{R}}D^{r} \sin x \big|_{x=0}
= \sin\!\left(\frac{\pi r}{2}\right).
\]
\textbf{Correction terms}. The integrand simplifies to
\[
\phi(r,x) = \sin\!\left(\frac{\pi r}{2}\right)\frac{x^{r}}{\Gamma(r+1)}.
\]
From Definition~\ref{def:correction_series}, the first three correction terms are:
\[
\begin{aligned}
\mathcal{E}_{\leq 2}[\sin](x)
={}&
\underbracket{0}_{\textstyle \mathcal{E}_0}
\;+\;
\underbracket{-\frac{\pi}{24}}_{\textstyle \mathcal{E}_1}
\;+\;
\underbracket{\frac{\pi}{1920}\!\left(4(\log x+\gamma)^{2}-\pi^{2}\right)}_{\textstyle \mathcal{E}_2}
\;+\;
\cdots.
\end{aligned}
\]

The vanishing of $\mathcal E_0$ reflects the common endpoint formula $\mathcal E_0=\frac12 f(0)$.

\begin{figure}[tbp]
  \centering
  \includegraphics[width=\linewidth]{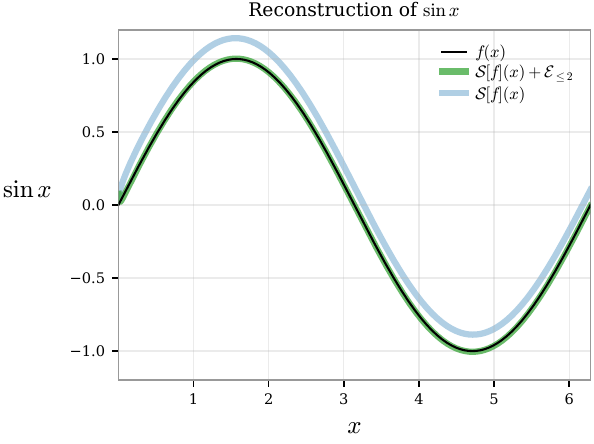}
  \caption{Spectral operator reconstruction of $\sin x$ showing the uncorrected operator and with Euler--Maclaurin corrections $\mathcal E_{\leq 2}$. Since $f(0)=0$, the common $\mathcal{E}_0$ correction term vanishes.}
  \label{fig:sine}
\end{figure}

\subsection{Rapidly decaying function: $f(x)=e^{-x^2}$}
\textbf{Fractional derivative}. For the Gaussian function, the Liouville (Fourier-multiplier) definition is admissible and yields the result below.
\[
{}^{\mathrm{L}}_{\mathbb{R}}D^{r} e^{-x^{2}}
=
\mathcal{F}^{-1}\!\left[(i\xi)^{r} e^{-\xi^{2}/4}\right](x),
\qquad\text{so that}\qquad
{}^{\mathrm{L}}_{\mathbb{R}}D^{r} e^{-x^{2}}\big|_{x=0}
=
\frac{\Gamma(r+1)}
{\Gamma\!\left(\frac{r}{2}+1\right)}
\cos\!\left(\frac{\pi r}{2}\right).
\]
\textbf{Correction terms}. The integrand simplifies to
\[
\phi(r,x)
=
\frac{\cos\!\left(\frac{\pi r}{2}\right)}{\Gamma\!\left(\frac{r}{2}+1\right)}\,x^{r}.
\]
We evaluated the required derivatives numerically using \texttt{mpmath.diff}. The numerical reconstruction behaves the same way as the exponential function. The shared $\mathcal{E}_0$ correction term provides the dominant improvement, while further corrections up to $\mathcal E_{\leq 2}$ refine the residual. Because the behavior parallels that of Figure~\ref{fig:exponential}, we omit the corresponding plots.

\subsection{Special function with oscillatory decay: $f(x)=J_0(x)$}
\textbf{Fractional derivative}. For the Bessel function, the Liouville (Fourier-multiplier) definition is admissible and yields the result below.
\[
{}^{\mathrm{L}}_{\mathbb{R}}D^{r} J_{0}(x)
=
\mathcal{F}^{-1}\!\left[(i\xi)^{r}
\frac{2}{\sqrt{1-\xi^{2}}}\mathbf{1}_{|\xi|<1}\right](x),
\qquad\text{so that}\qquad
{}^{\mathrm{L}}_{\mathbb{R}}D^{r} J_{0}(x)\big|_{x=0}
=
\frac{\Gamma(r+1)}
{2^{r}\Gamma\!\left(\frac{r}{2}+1\right)^{2}}
\cos\!\left(\frac{\pi r}{2}\right).
\]
\textbf{Correction terms}. The integrand simplifies to
\[
\phi(r,x)
=
\frac{\cos\!\left(\frac{\pi r}{2}\right)}{2^{r}\Gamma\!\left(\frac{r}{2}+1\right)^{2}}\,x^{r}.
\]
We evaluated the required derivatives numerically using \texttt{mpmath.diff}. Despite the oscillatory decay of $J_0$, the reconstruction behaves the same way as the previous examples. The spectral operator closely captures the global behavior, while Euler--Maclaurin corrections refine the remaining error. We omit representative plots because the behavior is qualitatively similar to Figures~\ref{fig:exponential} and \ref{fig:sine}.

\subsection{Function with a real pole: $f(x) = 1/(1-x)$}
This function has a Maclaurin radius of convergence $R=1$ due to the real pole at $x = 1$.

\textbf{Fractional derivative}. For the function with a real pole, we use the left-sided Liouville derivative based at $-\infty$, which gives
\[
\prescript{\mathrm{L}}{-\infty}{D}^{r}\!\left(\frac{1}{1-x}\right)
= \frac{\Gamma(r+1)}{(1-x)^{r+1}},
\qquad\text{so that}\qquad
\prescript{\mathrm{L}}{-\infty}{D}^{r}f(0) = \Gamma(r+1).
\]
\textbf{Correction terms}. The integrand simplifies to
\[
\phi(r,x) = x^{r}.
\]
From Definition~\ref{def:correction_series}, the first three correction terms are:
\[
\mathcal{E}_{\leq 2}\!\left[\tfrac{1}{1-x}\right](x)
=
\underbracket{\frac{1}{2}}_{\textstyle \mathcal{E}_0}
\;+\;
\underbracket{- \frac{1}{12}\,\log x}_{\textstyle \mathcal{E}_1}
\;+\;
\underbracket{\frac{1}{720}\,(\log x)^{3}}_{\textstyle \mathcal{E}_2}
\;+\;
\cdots.
\]
The spectral operator and the classical Maclaurin series share the same convergence boundary $R=1$ here, though the underlying mathematics differ.
\begin{figure}[tbp]
  \centering
  \includegraphics[width=\linewidth]{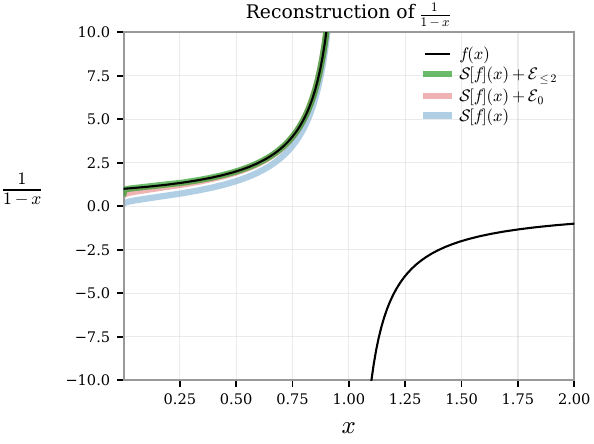}
  \caption{Spectral operator reconstruction of $1/(1-x)$ showing the uncorrected operator, with the common correction $\mathcal E_0$, and with Euler--Maclaurin corrections $\mathcal E_{\leq 2}$.}
  \label{fig:pole}
\end{figure}

\subsection{Function with complex poles: $f(x)=1/(1+x^{2})$}
The Lorentzian function pairs naturally with the real-pole example above: it has complex poles at $x=\pm i$ instead of a real one, so its Maclaurin series also has radius of convergence $R=1$.

\textbf{Fractional derivative}. For the Lorentzian function, the Liouville (Fourier-multiplier) definition is admissible and yields the result below.
\[
{}^{\mathrm{L}}_{\mathbb{R}}D^{r}f(x)
=
\mathcal{F}^{-1}
\!\left[
(i\xi)^r \pi e^{-|\xi|}
\right](x),
\qquad\text{so that}\qquad
{}^{\mathrm{L}}_{\mathbb{R}}D^{r}f(0)
=
\Gamma(r+1)\cos\!\left(\frac{\pi r}{2}\right).
\]
\textbf{Correction terms}. The integrand simplifies to
\[
\phi(r,x)
=
x^r\cos\!\left(\frac{\pi r}{2}\right).
\]
From Definition~\ref{def:correction_series}, the first three correction terms are:
\[
\mathcal{E}_{\leq 2}\!\left[\frac{1}{1+x^2}\right](x)
=
\underbracket{\frac12}_{\textstyle \mathcal{E}_0}
\;+\;
\underbracket{-\frac{1}{12}\log x}_{\textstyle \mathcal{E}_1}
\;+\;
\underbracket{
\frac{1}{720}
\left(
\log^3 x
-
\frac{3\pi^2}{4}\log x
\right)
}_{\textstyle \mathcal{E}_2}
\;+\;\cdots.
\]
As with the real pole example, both the spectral operator and the classical Maclaurin representations share the same convergence boundary $R=1$, though the underlying mathematics differ.

\begin{figure}[tbp]
  \centering
  \includegraphics[width=\linewidth]{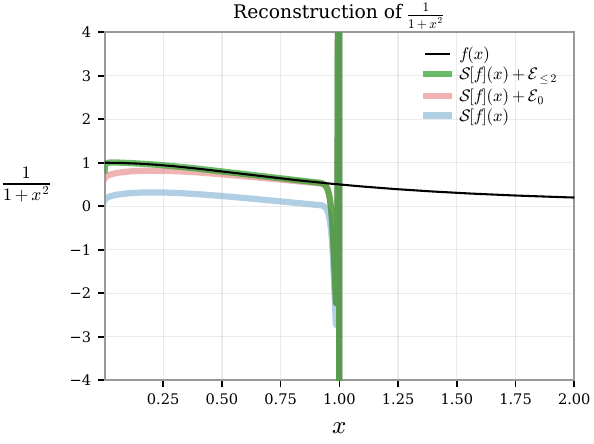}
  \caption{Spectral operator reconstruction of the Lorentzian function $1/(1+x^2)$ showing the uncorrected operator, with the common correction $\mathcal E_0$, and with Euler--Maclaurin corrections $\mathcal E_{\leq 2}$.}
  \label{fig:lorentzian}
\end{figure}

\begin{table}[tbp]
  \centering
  \caption{Mean absolute reconstruction errors for the uncorrected spectral operator, two Euler--Maclaurin correction levels, and the Abel--Plana correction integral. Errors were computed on a uniform grid over $0<x\leq 0.9$, emphasizing the region near $x=0$ where the sum--integral discrepancy is largest.}
  \label{tab:mae}
\begin{tabular}{lcccc}
\toprule
Function
& $\mathcal S$
& $\mathcal S+\mathcal E_0$
& $\mathcal S+\mathcal E_{\leq 2}$
& Abel--Plana \\
\midrule
$e^x$
& $0.546$& $0.0602$& $0.00142$& $\leq10^{-99}$\\

$\sin x$
& $0.139$& $0.139$& $0.00456$& $\leq10^{-99}$\\

$e^{-x^2}$
& $0.571$& $0.0743$& $0.00181$& $\leq10^{-99}$\\

$J_0(x)$
& $0.609$& $0.109$& $0.00337$& $\leq10^{-99}$\\

$1/(1-x)$
& $0.585$& $0.0851$& $0.00212$& $\leq10^{-99}$\\

$1/(1+x^2)$
& $0.594$& $0.0941$& $0.00181$& $\leq10^{-99}$\\
\bottomrule
\end{tabular}
\par\vspace{2mm}
\footnotesize\textit{Note.} Abel--Plana residuals were computed at 100-digit working precision against the correction integral of Proposition~\ref{prop:abel_plana_discrepancy}; all values fell at or below the working precision, consistent with the exact identity of that proposition.
\end{table}

The common correction $\mathcal E_0$ produces the largest reduction in error across the examples, except for $\sin x$, where $\mathcal E_0$ vanishes. Additional Euler--Maclaurin terms provide further, smaller reductions.

\subsection{Monomials}
\label{sec:monomials}
We now examine a degenerate case that tests the boundary of the admissibility conditions. For $f(x)=x^k$, the classical derivative data are concentrated at a single nonzero order,
\[
f^{(n)}(0)
=
\begin{cases}
0, & n \neq k,\\
k!, & n = k,
\end{cases}
\]
and so its Maclaurin series collapses to just one term. Similarly, the Riemann--Liouville derivative gives
\[
D^r x^k
=
\frac{\Gamma(k+1)}
{\Gamma(k-r+1)}
x^{k-r},
\qquad\text{so that}\qquad
D^r x^k(0)
=
\begin{cases}
0, & 0\le r<k,\\
k!, & r=k,\\
0, & r\in\{k+1,k+2,\ldots\},\\
\pm\infty, & r>k,\quad r\notin\mathbb N.
\end{cases}
\]
This appears to concentrate the nonzero derivative-order data at a single order, which has zero Lebesgue measure. However, for $r>k$, the Riemann--Liouville derivative at the origin fails the finiteness condition, so the spectral reconstruction operator is not defined for monomials under this choice of derivative. This is exactly the exclusion the finiteness condition of Section~\ref{sub:frac_deriv_admissibility} anticipates. In this case, the admissibility conditions actively rule out the one case where fractional derivative data fail to remain finite.

\section{Conclusions}

This paper introduces the spectral reconstruction operator $\mathcal{S}[f]$ as an independent mathematical object and establishes its foundational properties. We proved that $\mathcal{S}[e^x](x) \sim e^x$ as $x \to \infty$, the first proof that the operator asymptotically reconstructs its target function. We further showed that the Maclaurin summand and the spectral integrand are evaluations of the same underlying function of derivative order, making their discrepancy a sum--integral discrepancy. Classical summation formulas quantify this discrepancy under their respective hypotheses; we use the Euler--Maclaurin and Abel--Plana formulas here as representatives, sharing the same leading correction $\mathcal{E}_0=\frac{1}{2}f(0)$, which accounts for the dominant portion of the discrepancy. Under stronger analyticity conditions, Abel--Plana gives an exact integral representation with no remainder term. We verify numerically that these hypotheses hold in every example considered here. In particular, the Abel--Plana representation matches $f$ to the full precision of our computations throughout, as reported in Table~\ref{tab:mae}.

The classical Maclaurin series draws only on local information: integer derivatives at the origin. The derivatives entering $\mathcal{S}[f]$, by contrast, are non-local: the left-sided Liouville derivative draws on function values over the half-line extending to $-\infty$, while the Liouville (Fourier-multiplier) derivative draws on function values over all of $\mathbb R$. Either way, an infinite stretch of information about $f$ collapses into the single anchored value $D^r f(0)$. This combination of non-local derivative data and continuous integration over all nonnegative real orders distinguishes $\mathcal{S}[f]$ from all prior discrete fractional expansions.

The examples confirm stable reconstruction across entire, oscillatory, rapidly decaying, special-function, and finite-radius function classes. The uncorrected operator already provides useful approximations where $f(x)$ is large, since in the examples studied, the sum--integral discrepancy is bounded and decays with $x$. Summation corrections systematically account for the discrepancy concentrated near the origin, as shown in Table~\ref{tab:mae}. In the finite-radius examples studied (one with a real pole, one with a complex-conjugate pole pair), the operator encounters the same convergence boundary as the Maclaurin series, though through different mathematics: the Maclaurin partial sums diverge as a discrete series, while the spectral integrand fails to converge beyond this boundary.

The continuous-order data $D^r f(0)$ are not canonical: different admissible fractional derivative families produce different order spectra for the same function, and the reconstruction might depend on this choice. The admissibility conditions introduced here (recovery of the integer derivative ladder, the semigroup property, and finiteness of $D^r f(0)$ for all $r>0$) constrain but do not uniquely determine the derivative family. Which families are natural for a given function class, and whether the operator's behavior is stable across admissible choices, remain open questions. We have not fully characterized the conditions under which the spectral integral converges and produces meaningful reconstruction, beyond the case of $e^x$ proved here.

Several questions merit further investigation. The asymptotic reconstruction of $e^x$ relies on Laplace's method applied to the spectral integrand; extending this argument to other function classes would strengthen the theoretical foundation. We do not know which admissible families support the complex analyticity conditions required for exact integral representations.

\section*{Funding}
This research did not receive any specific grant from funding agencies in the public, commercial, or not-for-profit sectors.

\section*{Declaration of competing interest}
The author declares that he has no known competing financial interests or personal relationships that could have appeared to influence the work reported in this paper.

\section*{Data availability}
The code used to perform the numerical computations and generate the figures in this paper is not yet publicly available. It will be deposited in a public repository upon acceptance of this manuscript, and is available from the author upon reasonable request in the meantime.

\section*{Declaration of Generative AI and AI-assisted technologies in the writing process}
During the preparation of this work, the author used Claude (Anthropic) to improve the language, clarity, and structure of the manuscript text. After using this tool, the author reviewed and edited the content as needed and takes full responsibility for the content of the publication.

\bibliographystyle{plain}
\bibliography{references}

@article{Alquran2023,
  author  = {Alquran, M.},   
  title   = {The amazing fractional {Maclaurin} series for solving different types of fractional mathematical problems that arise in physics and engineering},
  journal = {Part. Differ. Equ. in Appl. Math.},
  volume  = {7},
  pages   = {100506},
  year    = {2023},
  doi     = {10.1016/j.padiff.2023.100506}
}

@article{FernandezBaleanu2018,
  author  = {Fernandez, A. and Baleanu, D.},
  title   = {The mean value theorem and {Taylor}'s theorem for fractional derivatives with {Mittag--Leffler} kernel},
  journal = {Adv. in Differ. Equ.},
  volume  = {2018},
  pages   = {86},
  year    = {2018},
  doi     = {10.1186/s13662-018-1543-9}
}

@misc{Gladkina2017,
  author  = {Gladkina, A. and Shchedrin, G. and {Al Khawaja}, U. and Carr, L. D.},
  title   = {Expansion of fractional derivatives in terms of an integer derivative series: Physical and numerical applications},
  year    = {2017},
  note    = {Preprint, Version 2. arXiv:1710.06297.},
  eprint  = {1710.06297},
  archivePrefix = {arXiv},
  doi     = {10.48550/arXiv.1710.06297}
}

@book{NIST:DLMF,
  title        = {NIST Digital Library of Mathematical Functions},
  editor       = {Olver, F. W. J. and {Olde Daalhuis}, A. B. and Lozier, D. W. and
                  Schneider, B. I. and Boisvert, R. F. and Clark, C. W. and
                  Miller, B. R.},
  year         = {2023},
  publisher    = {National Institute of Standards and Technology},
  url          = {https://dlmf.nist.gov/}
}

@article{Osler1971,
  author  = {Osler, T. J.},
  title   = {Taylor's series generalized for fractional derivatives and applications},
  journal = {SIAM J. Math. Anal.},
  volume  = {2},
  number  = {1},
  pages   = {37--48},
  year    = {1971},
  doi     = {10.1137/0502004}
}

@article{Ross1977,
  author  = {Ross, B.},
  title   = {The Development of Fractional Calculus 1695--1900},
  journal = {Hist. Math.},
  volume  = {4},
  number  = {1},
  pages   = {75--89},
  year    = {1977},
  doi     = {10.1016/0315-0860(77)90039-0}
}

@misc{Wei2020,
  author        = {Wei, Y. and Chen, Y. and Gao, Q. and Wang, Y.},
  title         = {Infinite series representation of fractional calculus: theory and applications},
  year          = {2020},
  note          = {Preprint, Version 3. arXiv:1901.11134.},
  eprint        = {1901.11134},
  archivePrefix = {arXiv},
  doi           = {10.48550/arXiv.1901.11134}
}

@article{Weyl1917,
  author  = {Weyl, H.},
  title   = {Bemerkungen zum {Begriff} des {Differentialquotienten} 
             gebrochener {Ordnung}},
  journal = {Vierteljahresschrift der Naturforschenden Gesellschaft 
             in {Z}{\"u}rich},
  volume  = {62},
  pages   = {296--302},
  year    = {1917}
}

@article{ZitaneTorres2023,
  author  = {Zitane, H. and Torres, D. F. M.},
  title   = {Generalized {Taylor}'s formula for power fractional derivatives},
  journal = {Bol. Soc. Mat. Mex.},
  volume  = {29},
  number  = {3},
  pages   = {Paper No. 68, 14pp},
  year    = {2023},
  doi     = {10.1007/s40590-023-00540-0}
}

\end{document}